\documentclass[12pt,reqno]{article}

\usepackage[usenames]{color}
\usepackage{amssymb}
\usepackage{graphicx}
\usepackage{amscd}

\usepackage[colorlinks=true,
linkcolor=webgreen,
filecolor=webbrown,
citecolor=webgreen]{hyperref}

\definecolor{webgreen}{rgb}{0,.5,0}
\definecolor{webbrown}{rgb}{.6,0,0}

\usepackage{color}
\usepackage{fullpage}
\usepackage{float}

\usepackage{graphics,amsmath,amssymb}
\usepackage{amsthm}
\usepackage{amsfonts}
\usepackage{latexsym}
\usepackage{epsf}

\makeatletter

\theoremstyle{plain}
\newtheorem{theorem}{Theorem}
\newtheorem{lemma}[theorem]{Lemma}

\theoremstyle{definition}
\newtheorem{definition}[theorem]{Definition}

\newcommand{\summ}{\sum\limits}

\theoremstyle{definition}

\makeatother
\newcommand{\fxy}[3]{f_{#1}({#2},{#3})}
\newcommand{\f}[2]{f_{#1}({#2})}  
\begin{document}

\begin{center}
\uppercase{\bf Enumerating triangulations by parallel diagonals}
\vskip 20pt
{\bf Alon Regev}\\
{Department of Mathematical Sciences, Northern Illinois University, DeKalb, Illinois}\\
{\tt regev@math.niu.edu}\\

\end{center}


%

\thispagestyle{empty}
\baselineskip=12.875pt
\vskip 30pt

\section{Introduction}

The purpose of this note is to enumerate triangulations of a regular convex polygon according to the number of diagonals parallel to a fixed edge.
This enumeration is of interest because it provides insight into the ``shape of a typical triangulation" and because of its connection to the Shapiro convolution.

We consider a triangulation of an $n$-gon as a labeled graph with vertices $0,1,\ldots, n-1$\footnote{For convenience, the vertex $0$ of an $n$-gon is sometimes also labeled $n$.} and edges denoted $xy$ for distinct vertices $x$ and $y$. The edges include $n$ sides $01, 12, \ldots ,(n-1)0$ and $n-3$ diagonals.

\begin{figure}[h]
\begin{center}
\epsfxsize=6.5in
\leavevmode\epsffile{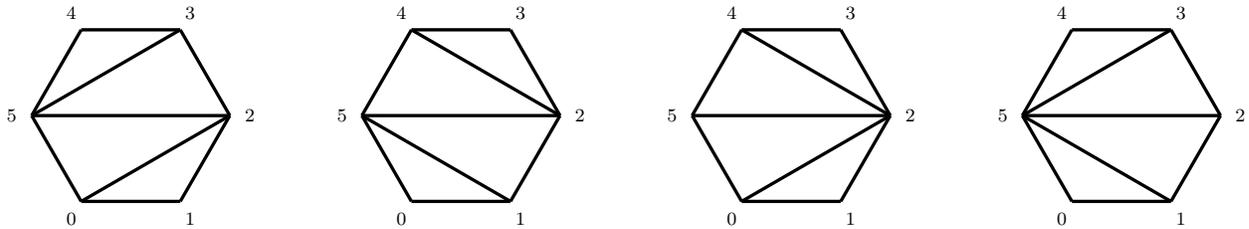}
\end{center}
\caption{The triangulations of a hexagon that include one diagonal parallel to $01$.}
\label{example1}
\end{figure}

\begin{definition}
Let $\fxy{xy}{n}{k}$ be the number of triangulations of a regular $n$-gon which include exactly $k$ diagonals parallel to the edge $xy$. Also denote $\fxy{xy}{n}{0}$ by $\f{xy}{n}$.
\end{definition}
For example, there are $14$ triangulations of a hexagon, $4$ of which include a diagonal parallel to $01$ (see Figure \ref{example1}). The remaining $10$ triangulations all have zero diagonals parallel to $01$. Therefore $\fxy{01}{6}{1}=4$ and $\f{01}{6}=10$.

Given $n$ and $k$, by symmetry $\fxy{xy}{n}{k}$ depends only on the value of $y-x$ modulo $n$ and not on the specific choice of $x$ and $y$. Furthermore, two edges $ab$ and $cd$ in a triangulation of an $n$-gon are parallel if and only if $a+b$ and $c+d$ are congruent modulo $n$.
It follows that for all $n$, $x$ and $y$,

\begin{equation*}
\fxy{xy}{2n}{k}=
\begin{cases}
\fxy{01}{2n}{k}, \qquad \text{ if $x+y$ is odd,}\\
\fxy{02}{2n}{k}, \qquad \text{ if $x+y$ is even,}\\
\end{cases}
\end{equation*}
and
\begin{equation*}
\fxy{xy}{2n+1}{k}=\fxy{01}{2n+1}{k}.
\end{equation*}

The question at hand is thus reduced to finding $\fxy{01}{2n}{k}$, $\fxy{02}{2n}{k}$ and \\$\fxy{01}{2n+1}{k}$.
Theorems \ref{mainthm} and \ref{k>0} below provide explicit formulas for these functions when $k=0$ and when $k>0$, respectively. These formulas are given in terms of the Catalan numbers
\[C_n={1\over n+1} {2n\choose n}.\]
Recall that there are $C_{n-2}$ triangulations of an $n$-gon. 
Therefore for all $n$, $x$ and $y$,
\begin{align}\label{rec1}
\summ_{k\ge 0}\fxy{xy}{n}{k}=C_{n-2}.
\end{align}
The recursion relation
\begin{align}\label{standard}
\summ_{i=0}^nC_iC_{n-i}=C_{n+1}
\end{align}
implies the identity
\begin{align}\label{hcon}
\summ_{i=0}^n C_{2i}C_{2n+1-2i}={1\over 2} C_{2n+2},
\end{align}
which is used below.
We also make use of the Shapiro convolution identity:
\begin{align}\label{shap}
\summ_{j=0}^nC_{2j}C_{2n-2j}=4^nC_n.
\end{align}
Andrews \cite{A} recently gave several proofs of \eqref{shap} and its $q$-analogs, with one of these proofs being purely combinatorial (however, finding a simple bijective proof of \eqref{shap} is still an open problem).

\section{Avoiding diagonals of a fixed direction}

We begin by enumerating the triangulations that avoid all diagonals parallel to a fixed edge.
\begin{theorem}\label{mainthm}
For any $n\ge 2$,
\begin{align}\label{even01eq}
\f{01}{2n}=2C_{2n-3}
\end{align}
and
\begin{align}\label{even02eq}
\f{02}{2n}=C_{2n-1}+2C_{2n-2}-2^{2n-1}C_{n-1}.
\end{align}
For any $n\ge 1$,
\begin{align}\label{oddeq}
\f{01}{2n+1}=2^{2n-1}C_{n-1}-C_{2n-1}.
\end{align}
\end{theorem}
Equations \eqref{even01eq}, \eqref{even02eq}  and \eqref{oddeq} can be proved by induction on $n$; the base cases are easily verified.
\begin{figure}
\begin{center}
\epsfxsize=2.5in
\leavevmode\epsffile{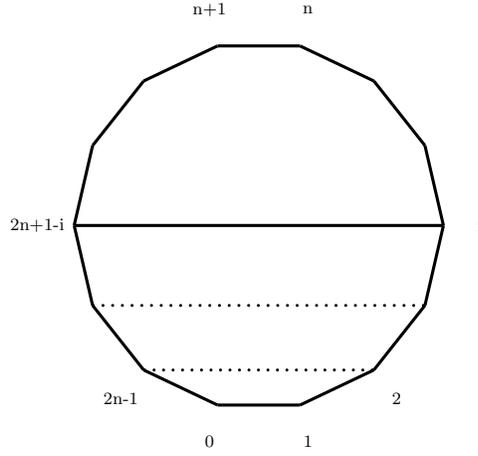}
\end{center}
\caption{Illustration of the proof of \eqref{even01eq}. The dotted lines represent the avoided diagonals.}\label{FigA}
\end{figure}

\subsection*{Proof of \eqref{even01eq}:}

Enumerate the triangulations of a $2n$-gon that include at least one diagonal parallel to $01$ according to the minimal number $i$, with $2\le i \le n-1$, such that $i(2n+1-i)$ is an edge of the triangulation
(see Figure \ref{FigA}).
The $(2n-2i+2)$-gon with vertices $i,i+1,\ldots, 2n+1-i$ can be triangulated in $C_{2n-2i}$ ways.
By induction the $2i$-gon with vertices $0,1,\ldots, i, 2n+1-i, 2n+2-i,\ldots ,2n-1$ can be triangulated in $2C_{ 2i-3}$ ways.
Subtracting these from the total number triangulations of a $2n$-gon gives
\begin{align*}
\f{01}{2n}&=C_{2n-2}-\summ_{i=2}^{n-1}C_{2n-2i}\cdot 2C_{2i-3}\\
&=C_{2n-2}-2\summ_{i=2}^nC_{2n-2i} C_{2i-3}+2C_0C_{2n-3}\\
&=2C_{2n-3},\\
\end{align*}
where in the last equality we have used \eqref{hcon}.
Another proof of \eqref{even01eq}, using a result of David Callan on Dyck paths, is outlined in Section \ref{remarks-sec}.

\begin{figure}
\begin{center}
\epsfxsize=2.5in
\leavevmode\epsffile{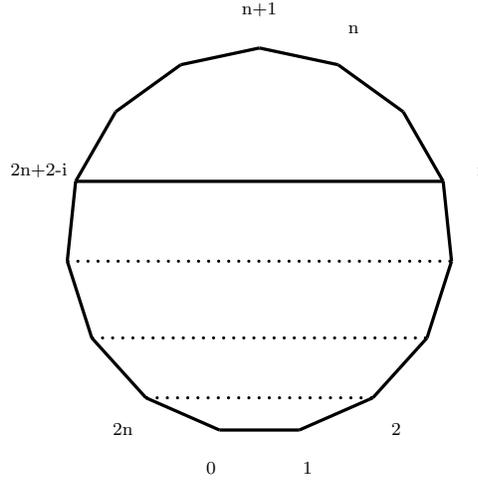}
\end{center}
\caption{Illustration of the proof of \eqref{oddeq}.}\label{FigB}
\end{figure}

\subsection*{Proof of \eqref{oddeq}:}

Enumerate the triangulations of a $(2n+1)$-gon that include at least one diagonal parallel to $01$ according to their diagonal $i(2n+2-i)$ with minimal $i$
(see Figure \ref{FigB}).
The $(2n-2i+3)$-gon with vertices $i, i+1, \ldots , 2n+2-i$ can be triangulated in $C_{2n-2i+1}$ ways.
By \eqref{even01eq}, the $2i$-gon with vertices $0,1, \ldots, i, 2n+2-i, 2n+3-i, \ldots , 2n$ can be triangulated in $2C_{2i-3}$ ways.
Therefore
\begin{align*}
\f{01}{2n+1}&=C_{2n-1}-\summ_{i=2}^nC_{2n-2i+1}\cdot 2C_{2i-3}\\
&=\summ_{j=0}^{2n-2}C_{j}C_{2n-2-j}-2\summ_{i=2}^nC_{2n-2i+1}C_{2i-3}\\
&=\summ_{j=0}^{2n-2}(-1)^jC_jC_{2n-2-j}.\\
\end{align*}

Thus by \eqref{standard} and \eqref{shap},
\begin{align*}
\f{01}{2n+1}&=2\summ_{j=0}^{n-1} C_{2j}C_{2n-2-2j}-\summ_{j=0}^{2n-2} C_jC_{2n-2-j}=2^{2n-1}C_{n-1}-C_{2n-1}.\\
\end{align*}

\begin{figure}
\begin{center}
\epsfxsize=2.5in
\leavevmode\epsffile{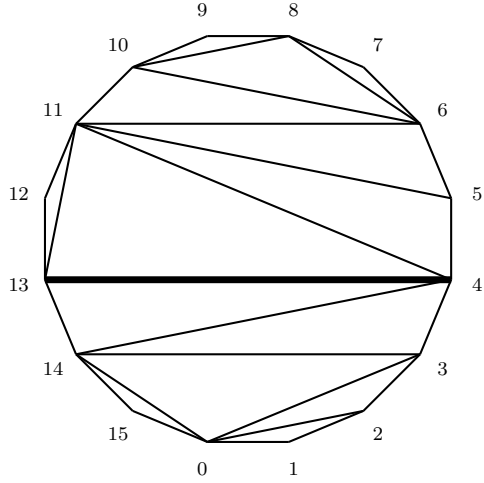}
\end{center}
\caption{Example illustrating the proof of Lemma \ref{lem}.}\label{PLF}
\end{figure}

The following lemma will be used in the proof of \eqref{even02eq}.
\begin{lemma}\label{lem}
For any $n\ge 2$,
\begin{equation}\label{lemeq}
\summ_{i=1}^{n-1} 2^{2i-1}C_{i-1}C_{2n-1-2i}=4^{n-1}C_{n-1}-C_{2n-2}
\end{equation}
\begin{proof}
Let $h(n)$ be the number of triangulations of a $2n$-gon together with a marking either on one of the sides $01$, $n(n+1)$ or on one of the diagonals $k(2n+1-k)$, with $2\le k\le n-1$, if any such diagonals are present. For example, Figure \ref{PLF} shows a triangulation of a $16$-gon with the diagonal $4(13)$ marked.
Consider the following two ways to enumerate these marked triangulations.

\begin{enumerate}
\item First mark the edge $(j+1)(2n-j)$, with $0 \le j\le n-1$. For example, the marked triangulation in Figure \ref{PLF} corresponds to $n=8$ and $j=3$.
 Then choose one of the $C_{2j}$ triangulations of the $(2j+2)$-gon with vertices $0,1,\ldots, j+1, 2n-j, 2n-j+1, \ldots ,2n-1$ and one of the $C_{2n-2j-2}$ triangulations of the $(2n-2j)$-gon with vertices $j+1,j+2,\ldots, 2n-j$. Thus there are  $C_{2j-2}C_{2(n-1)-2j}$ such marked  triangulations for each $j$. By \eqref{shap},
\begin{align}\label{shap2}
h(n)=\summ_{j=0}^{n-1}C_{2j}C_{2(n-1)-2j} = 4^{n-1}C_{n-1}.
\end{align}
\item
There are $C_{2n-2}$ marked triangulations whose edge $n(n+1)$ is the one marked.
The remaining marked triangulations can be enumerated according to the maximal $i$, with $1 \le i\le n-1$, such that $i(2n+1-i)$ is one of the diagonals in the triangulation (where the case $i=1$ corresponds to triangulations avoiding all diagonals parallel to $01$.) For example, in Figure \ref{PLF} we have $n=8$ and $i=6$.
For each such $i$, there are $h(i)$ marked triangulations of the $2i$-gon with vertices $0,1, \ldots, i, 2n-2i+1, \ldots, 2n-1$, and there are
 $\f{01}{2n-2i+2}$ triangulations of the $(2n-2i+2)$-gon with vertices $i,i+1,\ldots, 2n+1-i$ which avoid the diagonals $(i+1)(2n-i),\ldots , n(n+1)$.
 Thus by \eqref{even01eq} and \eqref{shap2},
\begin{align*}
h(n)&=C_{2n-2}+\summ_{i=1}^{n-1}h(i)\f{01}{2n-2i+2}\\
&=C_{2n-2}+\summ_{i=1}^{n-1}4^{i-1}C_{i-1}\cdot 2C_{2n-2i-1}.
\end{align*}
\end{enumerate}
Comparing this with \eqref{shap2} completes the proof.
\end{proof}
\end{lemma}

\subsection*{Proof of \eqref{even02eq}:}

\begin{figure}
\begin{center}
\epsfxsize=2.5in
\leavevmode\epsffile{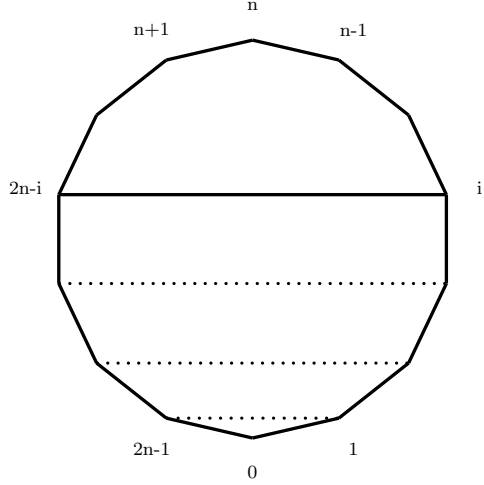}
\end{center}
\caption{Illustration of the proof of \eqref{even02eq}.}\label{FigC}
\end{figure}

For convenience we calculate $\f{1(2n-1)}{2n}=\f{02}{2n}$. Enumerate the triangulations of a $2n$-gon that include at least one diagonal parallel to $1(2n-1)$ according to their diagonal $i(2n-i)$ with minimal $i$, where $1\le i\le n-1$ (see Figure \ref{FigC}).
By \eqref{oddeq}, the $(2i+1)$-gon with with vertices $0, 1, \ldots,  i, 2n-i, 2n+1-i, \ldots,  2n-1$ can be triangulated in $2^{2i-1}C_{i-1}-C_{2i-1}$ ways. The $(2n-2i+1)$-gon with vertices $i,i+1,\ldots , 2n+2-i$ can be triangulated in  $C_{2n-2i-1}$ ways.
Therefore

\begin{align*}
\f{02}{2n}
&=C_{2n-2}-\summ_{i=1}^{n-1}(2^{2i-1}C_{i-1}-C_{2i-1})C_{2n-2i-1}\\
&=C_{2n-2}-\summ_{i=1}^{n-1}2^{2i-1}C_{i-1}C_{2n-1-2i}+\summ_{i=1}^{n-1}C_{2i-1}C_{2n-2i-1}\\
&=C_{2n-2}-(4^{n-1}C_{n-1}-C_{2n-2})+(C_{2n-1}-4^{n-1}C_{n-1})\\
&=2C_{2n-2}+C_{2n-1}-2^{2n-1}C_{n-1},\\
\end{align*}
where in the penultimate equality we have used \eqref{shap} and \eqref{lemeq}.

\section{Including a number of diagonals of a fixed direction}

The next theorem enumerates the triangulations with a fixed positive number of diagonals parallel to a fixed edge.
\begin{theorem}\label{k>0}
Let $n\ge 2$ and $k\ge 1$. Then
\begin{align}\label{k>0-even01}
\fxy{01}{2n}{k}= \summ_{\substack{i_1+\ldots +i_{k+1}=n-1}} 2^{k+1} C_{2i_1-1} C_{2i_2-1}\cdots C_{2i_{k+1}-1},
\end{align}
and
\begin{align}\label{k>0-even02}
\fxy{02}{2n}{k} =&\nonumber \\
\summ_{i_1+\ldots +i_{k+1}=n-1}&2^{k-1} (2^{i_1-1}C_{i_1-1}-C_{2i_1-1})(2^{i_2-1}C_{i_2-1}-C_{2i_2-1})C_{2i_3-1}C_{2i_4-1}\cdots C_{2i_{k+1}-1}. \end{align}
If $n,k\ge 1$ then
\begin{align}\label{k>0-odd}
\fxy{01}{2n+1}{k}=\summ_{i_1+\ldots +i_{k+1}=n} (2^{i_1-1}C_{i_1-1}-C_{2i_1-1})C_{2i_2-1}C_{2i_3-1} \cdots C_{2i_{k+1}-1}.
\end{align}
\begin{proof}
Consider a triangulation of an $2n$-gon which includes exactly $k$ diagonals parallel to $01$. These $k$ diagonals partition the $2n$-gon into $k+1$ triangulated polygons, and partition the $n-1$ edges $12, 23,\ldots, (n-1)n$ into $k+1$ corresponding parts consisting of $i_1,\ldots, i_{k+1}\ge 1$ edges. The number of vertices in each resulting polygon is $2i_j+2$ for all $j$, and each such polygon is triangulated with diagonals which are not parallel to one of its sides. Thus
\begin{align}\label{even01}
\fxy{01}{2n}{k}= \summ_{\substack{i_1+\ldots +i_{k+1}=n-1  \\ i_j\ge 1}}\f{01}{2i_1+2}\f{01}{2i_2+2}\cdots \f{01}{2i_{k+1}+2},
\end{align}
which together with \eqref{even01eq} proves equation \eqref{k>0-even01}.
By similar considerations,
\begin{align*}
\fxy{02}{2n}{k} =& \nonumber \\
\summ_{\substack{i_1+\ldots +i_{k+1}=n-1  \\ i_j\ge 1}}& \f{01}{2i_1+1}\f{01}{2i_2+1}
\f{01}{2i_3+2}\f{01}{2i_4+2}\cdots \f{01}{2i_{k+1}+2}, \qquad \qquad
\end{align*}
and
\begin{align*}
\fxy{01}{2n+1}{k}= \summ_{\substack{i_1+\ldots +i_{k+1}=n  \\ i_j\ge 1}}\f{01}{2i_1+1} \f{01}{2i_2+2}\f{01}{2i_3+2}\cdots \f{01}{2i_{k+1}+2}.
\end{align*}
The differences between these equations and \eqref{even01} result from considering the regions of the polygon which contain the vertices $0$ and $n$.
The proofs of \eqref{k>0-even02} and \eqref{k>0-odd} now follow from \eqref{even02eq} and \eqref{oddeq}, respectively.
\end{proof}
\end{theorem}

Note that a consequence of \eqref{rec1}, \eqref{even01eq} and \eqref{k>0-even01} is the Catalan identity
\begin{align}\label{catid1}
\summ_{\substack{k\ge 0\\ i_1+\ldots +i_{k+1}=n}} 2^{k+1} C_{2i_1-1} C_{2i_2-1}\cdots C_{2i_{k+1}-1}=C_{2n}.
\end{align}

Another Catalan identity can be obtained by considering the set of marked triangulations of a $2n$-gon described in the proof of Lemma \ref{lem}. If $k$ is the number of diagonals parallel to $01$ in a triangulation of a $2n$-gon, this triangulation corresponds to $k+2$ such marked triangulations. Thus by \eqref{k>0-even01},

\begin{align}\label{catid2}
\summ_{\substack{0 \le k \le n-2\\ i_1+\ldots +i_{k+1}=n-1}}(k+2) \,2^{k+1}C_{2i_1-1}C_{2i_2-1}\cdots C_{2i_{k+1}-1}=4^{n-1}C_{n-1}.
\end{align}

Combining \eqref{catid1} and \eqref{catid2} results in the identity

\begin{align*}
\summ_{\substack{1 \le k \le n-1\\ i_1+\ldots +i_{k+1}=n}}k\,2^kC_{2i_1-1}C_{2i_2-1}\cdots C_{2i_{k+1}-1}=2^{2n-1}C_n-C_{2n}.
\end{align*}

\section{Remarks}\label{remarks-sec}

The next theorem was proposed as a problem to the American Mathematical Monthly by David Callan in 2003, and a solution appeared in 2005.

\begin{theorem}\label{Cthm}\cite{C}
The number of Dyck $2n$-paths that avoid the points $(4k,0)$, $k=1,2,\ldots , n-1$ is twice the number of Dyck $(2n-1)$-paths.
\end{theorem}

Callan proved Theorem \ref{Cthm} using a bijection on Dyck paths. The result is equivalent to \eqref{even01eq}, since the Dyck paths in question are equinumerous with the triangulations of a $(2n-2)$-gon which avoid all diagonals parallel to $01$. To see this, compare the initial conditions for both sequences, and observe that the Dyck paths in question satisfy an analogous recursive relations to the ones given by equations \eqref{rec1} and \eqref{k>0-even01}.

Similarly, it can be shown that $\f{02}{2n}$ is equal to the number of Dyck $2n$-paths avoiding all points $(4k+2, 0)$ with $k=0,1,\ldots, n-1$, and that $\f{01}{2n+1}$ is equal to the number of Dyck $(2n+1)$-paths avoiding all points $(4k, 0)$  with $k=1,2,\ldots ,n-1$.

The relation with Dyck paths also gives another interpretation of these results in terms of triangulations. Using standard bijections between triangulations and Dyck paths, the points $(4k, 0)$ and $(4k+2,0)$ of a Dyck path correspond to the diagonals of the form $0(2k+1)$ and $0(2k)$, respectively, of a triangulation. This gives analogous results  to those of the present note, concerning the number of diagonals of this form instead of the number of diagonals parallel to a fixed edge.

The sequences $\f{02}{2n}$ and $\f{01}{2n+1}$ appear in \cite[A066357]{S} and \cite[A079489]{S}, respectively. The interpretation in terms Dyck paths is given there, along with other interpretations and several interesting properties. Callan (\cite{C-priv} and \cite[A066357]{S}) proved the analog of \eqref{even02eq} using generating functions. Barry \cite[A066357]{S}
gave an alternative formula for this sequence: 

\begin{equation}\label{Barry}
f_{02}(2n+2)={1\over n} \summ_{k=0}^n{4n \choose k} {3n-k-2 \choose n-k-1}.
\end{equation}

Callan used Dyck paths to prove that
\begin{equation*}
\f{02}{2n+2}=\summ_{k=1}^n \f{01}{2k+1}\f{01}{2(n-k)+1}.
\end{equation*}
Another relation between these sequences is evident from \eqref{even02eq} and \eqref{oddeq}:
\begin{equation}\label{2C}
\f{01}{2n+1}+\f{02}{2n}=2C_{2n-2}.
\end{equation}
A direct proof of \eqref{2C} may also be of interest.

\end{document}